\documentstyle[bezier,12pt]{article}

\input{emlines.sty}

\newcommand{\mn}{\sl}
\def\afterthmseparator{}
\makeatletter
\renewcommand{\@begintheorem}[2]{\trivlist
      \item[\hskip \labelsep{\bf #1\ #2\unskip\afterthmseparator}]\mn}
\renewcommand{\@opargbegintheorem}[3]{\trivlist
      \item[\hskip \labelsep{\bf #1\ #2\ (#3)\unskip\afterthmseparator}]\mn}
\makeatother
\newtheorem{theorem}{Theorem}[section]
\newtheorem{lemma}[theorem]{Lemma}

\newtheorem{proposition}[theorem]{Proposition}
\newtheorem{rem}[theorem]{Remark}
\newenvironment{remark}{\renewcommand{\mn}{\rm} \begin{rem}}{\end{rem}}
\newtheorem{probl}[theorem]{Problem}

\newtheorem{df}[theorem]{Definition}
\newenvironment{definition}{\renewcommand{\mn}{\rm} \begin{df}}{\end{df}}
\newtheorem{exmpl}[theorem]{Example}

\newcommand{\qed}{$\;\;\;\Box$}
\newenvironment{proof}{\par\smallbreak{\sl Proof.~}}
{\unskip\nobreak\hfill \qed \par\medbreak}

\newcommand{\function}[2]{:#1 \rightarrow #2}
\newcommand{\of}[1]{\left( #1 \right)}

\newcommand{\setdef}[2]{\left\{\hspace{0.5mm}#1:\hspace{0.5mm} #2\right\}}

\newcommand{\bbn}{ {\bf N} }
\newcommand{\bbr}{ {\bf R} }
\newcommand{\ga}{\alpha}
\newcommand{\gb}{\beta}

\newcommand{\gd}{\delta}
\newcommand{\eps}{\epsilon}
\newcommand{\Gg}{\Gamma}
\newcommand{\And}{\wedge}
\newcommand{\Or}{\vee}

\newcommand{\hide}[1]{}
\newcommand{\refeq}[1]{(\ref{#1})}

\newcommand{\EF}{Ehrenfeucht}
\newcommand{\hopl}{Hamiltonian outerplanar}
\newcommand{\hop}{HOP}
\newcommand{\edhop}{1-e.d.HOP}
\newcommand{\eedhop}{2-e.d.HOP}

\newcommand{\game}{\mbox{\sc Ehr}}
\newcommand{\D}[1]{\mbox{D} (#1)}
\newcommand{\DD}[2]{\mbox{D}_{#1} (#2)}
\newcommand{\cl}{{\cal C}}
\newcommand{\clo}{{\cal O}}
\newcommand{\Nest}{\mathop{\rm Nest}\nolimits}
\newcommand{\nest}[1]{\Nest(#1)}
\newcommand{\Alt}{\mathop{\rm alt}\nolimits}
\newcommand{\alt}[1]{\Alt(#1)}
\newcommand{\Sep}{\mathop{\it sep}\nolimits}
\newcommand{\sep}[1]{\Sep(#1)}
\newcommand{\Col}{\mathop{\rm color}\nolimits}
\newcommand{\col}[1]{\Col(#1)}
\newcommand{\flap}[1]{\langle #1 \rangle}

\title{
The First Order Definability\\ of Graphs with Separators\\
via the Ehrenfeucht Game
}

\author{Oleg Verbitsky\\
{\small Kyiv University, Ukraine}\\
{\small oleg@ov.litech.net}}

\date{30 March 2003}

\begin{document}
\maketitle

\begin{abstract}
We say that a first order formula $\Phi$ defines a graph $G$ if $\Phi$
is true on $G$ and false on every graph $G'$ non-isomorphic with $G$.
Let $\D G$ be the minimal quantifier rank of a such formula.
We prove that, if $G$ is a tree of bounded degree or a
Hamiltonian (equivalently, 2-connected) outerplanar graph,
then $\D G=O(\log n)$, where $n$ denotes the order of $G$.
This bound is optimal up to a constant factor.
If $h$ is a constant, for connected graphs with no minor $K_h$
and degree $O(\sqrt n/\log n)$, we prove the bound $\D G=O(\sqrt n)$.
This result applies to planar graphs and, more generally, to graphs
of bounded genus.

Our proof techniques are based on the characterization of the
quantifier rank as the length of the \EF\/ game on non-isomorphic graphs.
We use the separator theorems to design a winning strategy
for Spoiler in this game.
\end{abstract}

\section{Introduction}
{}We treat a graph $G$ as a structure with a single
anti-reflexive and symmetric binary predicate $E$ for the adjacency
relation of $G$. Every closed first order formula $\Phi$ with predicate
symbols $E$ and $=$ is either true or false on $G$.
The {\em quantifier rank\/} of $\Phi$ is the maximum
number of nested quantifiers in this formula (see Section \ref{s:prel} for
formal definitions).
Given two non-isomorphic graphs $G$ and $G'$, we say that $\Phi$
{\em distinguishes $G$ from $G'$} if $\Phi$ is true on $G$ but false on
$G'$. Let $\D{G,G'}$ denote the minimum quantifier rank of
a such formula.

The definability issues, studied in finite model theory,
are important for such areas in computer science as databases,
model checking, and descriptive complexity theory.
The number $\D{G,G'}$ should be considered the measure of
indistinguishability of the graphs by means of first order logic.
The larger $\D{G,G'}$ is, the harder is to find a difference between
$G$ and $G'$ expressible by a first order statement. A characterization
of $\D{G,G'}$ in terms of partial isomorphisms between $G$ and $G'$
and extensions thereof is given by Fra\"\i ss\'e~\cite{Fra}.
Equivalently in essence, Ehrenfeucht \cite{Ehr} characterizes $\D{G,G'}$
as the length of a Spoiler-Duplicator game on $G$ and $G'$,
which provides us with a nice and robust tool for estimation of $\D{G,G'}$.

One of the simplest examples of graphs whose distinguishion requires
large quantifier rank
is given by complete graphs. The lower bound $\D{K_n,K_{n+1}} > n$
is obvious from the Ehrenfeucht characterization and matches the
simple general upper bound $\D{G,G'}\le n+1$ for $G$ of order $n$ and $G'$
non-isomorphic with $G$. A bit less obvious, rather popular in the
literature, examples of the first order similarity are given by
paths and cycles. Two paths or two cycles of distinct length turn out
indistinguishable by formulas of logarithmic quantifier rank:
if $m>n$, then
\begin{equation}\label{eq:path}
\D{P_n,P_m}>\log_2(n-1)-2
\end{equation}
(e.g.\ \cite[Theorem 2.1.3]{Spe}) and
\begin{equation}\label{eq:cycle}
\D{C_n,C_m}>\log_2n
\end{equation}
(e.g.\ \cite[Proof of Theorem 2.4.2]{Spe} and \cite[Example 2.3.8]{EFl}).
These lower bounds are well  known to be tight up to an
additive constant (e.g.\ \cite[Theorem 2.1.2]{Spe}).

We say that a first order formula $\Phi$ {\em defines\/} a graph $G$ if
$\Phi$ distinguishes $G$ from all non-isomorphic graphs.
Let $\D G$ denote the minimum quantifier rank of a such formula.
Complementing \refeq{eq:path} and \refeq{eq:cycle}, it is not hard to
show that
\begin{equation}\label{eq:pc}
\D{P_n}<\log_2n+3
\mbox{\ \ and\ \ }
\D{C_n}<\log_2n+3
\end{equation}
(cf.\ Remark \ref{rem:path}).
We here generalize \refeq{eq:pc} both for paths and cycles by showing that
\begin{equation}\label{eq:log}
\D G=O(\log n)
\end{equation}
for $G$ of order $n$ being, respectively, a tree of bounded degree or
a \hopl\/ graph.

If we put no restriction on the degree, $\D G$ for a tree may be much
bigger. For example, for stars we have
\begin{equation}\label{eq:star}
\D{K_{1,n-1},K_{1,n}}=n.
\end{equation}
The same example shows that the hamiltonicity is essential for
the definability of outerplanar graphs with small (logarithmic)
quantifier rank.
Note, however, that for outerplanar graphs the hamiltonicity should not
be considered a strong structural restriction. It is well known that
every 2-connected component of an outerplanar graph is an outerplanar
graph (see e.g. \cite{Har}) and it is easy to see that
2-connected outerplanar graphs are precisely \hopl\/ graphs.

It is worth noting that the defining formulas assumed by \refeq{eq:log}
have a restricted logical structure. We say that a first order formula
$\Phi$ is in {\em negation normal form\/} if the connective $\neg$
occurs in $\Phi$ only in front of atomic subformulas. If $\Phi$ is
such a formula, its {\em alternation number\/} is the maximum number
of alternations of $\exists$ and $\forall$ in a sequence of nested
quantifiers of $\Phi$.
Our proof of \refeq{eq:log} produces defining formulas
in negation normal form whose alternation number is at most~2.

The proof of \refeq{eq:log} is based on the analysis of the \EF\/
game on non-isomorphic graphs $G$ and $G'$, where $G$ is a tree
of bounded degree or \hopl. The analysis is uniform for both classes
of graphs and uses only the existence of a small separator in $G$.
A {\em separator\/} of $G$ is a set of vertices whose removal
splits $G$ into connected components each having at most $\eps n$
vertices. Cai, F\"urer, and Immerman \cite{CFI} prove that,
if $G$ has separator of size $O(n^\delta)$, $0<\delta<1$,
then $G$ is definable by a formula with counting quantifiers
that has $O(n^\delta)$ variables. For an integer $m$,
a {\em counting quantifier\/} $\exists^{\ge m}x\Psi$ means that there
are at least $m$ vertices $x$ for which the statement $\Psi$ holds.
Note that we use separators in a more complicated situation than they
were used in \cite{CFI}, where the power of counting quantifiers was
essential.

To prove the logarithmic bound \refeq{eq:log}, it is essential that
trees and outerplanar graphs have separators of constant size.
For some larger classes, as planar graphs, graphs of bounded genus,
and, most generally, graphs with an excluded minor, there exist
separators of size $O(\sqrt n)$~\cite{LTa,GHT,AST}.
We are able to adapt our techniques for connected graphs in these classes
and prove the bound $\D G=O(\sqrt n)$ under the restriction of
the maximum vertex degree of $G$ to $O(\sqrt n/\log n)$.
The connectedness is here essential: for example,
$\D{Triv_{m,2m},Triv_{m-1,2m+2}}=2m$, where $Triv_{a,b}$ is a graph
with $a$ isolated edges and $b$ isolated vertices. The restriction
of the maximum vertex degree is essential by \refeq{eq:star}.
Note that the bounded degree, on its own right, does not bound
$\D G$ so much. Cai, F\"urer, and Immerman \cite{CFI} construct
a sequence of graphs $G$ whose maximum degree is 4 but nevertheless
$\D G=\Omega(n)$.

{\bf Related work.}
In \cite{PVV} we prove that, if $G$ and $G'$ are non-isomorphic graphs
of the same order $n$, then $\D{G,G'}\le(n+3)/2$. Simple examples show
that this bound is best possible up to an additive constant of 1.
Thus, trees of bounded degree and Hamiltonian outerplanar graphs
should be considered classes of graphs definable with low quantifier
rank. Another important characteristics is the minimum number of variables
used in a distinquishing or defining formula (different occurrences
of the same variable are not counted). Note that this number
does not exceed the minimum quantifier rank.
Immerman and Kozen \cite{IKo} prove that every tree of degree
$d$, represented by a child-parent relation between vertices,
is definable with at most $d$ variables.
Grohe \cite{Gro1,Gro2} shows that planar graphs
and, more generally, graphs of bounded genus are definable by formulas
with counting quantifiers using only constantly many variables.

The paper is organized as follows. Section \ref{s:prel} contains
the relevant definitions from graph theory and logic.
We present our approach in Section \ref{s:str} and apply it to trees
of bounded degree and \hopl\/ graphs in Section \ref{s:hop}.
In Section \ref{s:minor} we handle the classes of graphs of bounded
degree with separators.

\section{Preliminaries}\label{s:prel}

\subsection{Graphs}

Given a graph $G$, we denote its vertex set by $V(G)$.
Throughout the paper, unless stated otherwise, $n$ denotes
the order of $G$, that is, $n=|V(G)|$. Sometimes the order of $G$
will be denoted by $|G|$. The neighborhood of a vertex $v$ consists of
all vertices adjacent to $v$ and is denoted by $\Gg(v)$.
The neighborhood of a set $X\subseteq V(G)$ is defined by
$\Gg(X)=\bigcup_{v\in X}\Gg(v)\setminus X$.
The {\em degree\/} of a vertex $v$ is the number of vertices in $\Gg(v)$.
The maximum vertex degree of $G$ is denoted by~$\Delta(G)$.
We say that graphs in a class $\cl$ have bounded degree if
$\Delta(G)=O(1)$ for $G$ in $\cl$.
The distance between two vertices of a graph G, $u$ and $v$, is
equal to the shortest length of a path between $u$ and $v$
and denoted by $d(u,v)$. The distance from a vertex $u\in V(G)$
to a set $X\subseteq V(G)$ is defined by $d(u,X)=\min_{v\in X}d(u,v)$.

A {\em Hamilton cycle\/} in a graph $G$ is an ordering $v_1,\ldots,v_n$
of $V(G)$ such that $v_i$ and $v_{i+1}$ are adjacent for all $i<n$
and adjacent are also $v_n$ and $v_1$. A graph with a Hamilton cycle
is called {\em Hamiltonian}.

An {\em outerplanar\/} graph is a planar graph embeddable in plane
so that all vertices lie on the border of the same face.

If $X\subseteq V(G)$, then $G[X]$ denotes the subgraph induced by $G$ on $X$.
The result of removal of all vertices in $X$ from $G$ is denoted by
$G\setminus X$, that is, $G\setminus X=G[V(G)\setminus X]$.

\begin{definition}
{\bf \cite{AST}}
If $X\subseteq V(G)$, we call a connected component of $G\setminus X$
an {\em $X$-flap}.
\end{definition}

We write $G\cong G'$ if graphs $G$ and $G'$ are isomorphic.
A one-to-one map $\phi\function{X}{X'}$, where $X\subseteq V(G)$
and $X'\subseteq V(G')$, is a {\em partial isomorphism\/} from $G$
to $G'$ if $\phi$ is an isomorphism from $G[X]$ to $G'[X']$.

A {\em colored graph\/} is a structure that, in addition to the
anti-reflexive and symmetric binary relation, has countably many
unary relations $C_i$, $i\ge 1$. The truth of $C_i(v)$ for a vertex $v$
is interpreted as coloration of $v$ in color $i$.
We consider finite colored graphs,
whose vertices can have only finitely many colors.
Most graph-theoretic notions carry over colored graphs literally.
For example, a colored graph is connected iff so is its underlying
graph. Note only that
an isomorphism of colored graphs preserves the adjacency relation
and, moreover, matches a vertex of one graph to an equally colored vertex
of the other graph. If $\cl$ is a class of graphs, $\col\cl$ will
denote the class of colored graphs whose underlying graphs are in~$\cl$.

\subsection{Logic}

First order formulas are assumed to be over the set of connectives
$\{\neg,\And,\Or\}$.

\begin{definition}
A {\em sequence of quantifiers\/} is a finite word over the alphabet
$\{\exists,\forall\}$. If $S$ is a set of such sequences, then
$\exists S$ (resp.\ $\forall S$) means the set of concatenations
$\exists s$ (resp.\ $\forall s$) for all $s\in S$. If $s$ is a sequence
of quantifiers, then $\bar s$ denotes the result of replacement of all
occurrences of $\exists$ to $\forall$ and vice versa in $s$. The set $\bar S$
consists of all $\bar s$ for $s\in S$.

Given a first order formula $\Phi$, its set of {\em sequences of nested
quantifiers\/} is denoted by $\nest\Phi$ and defined by induction as
follows:
\begin{enumerate}
\item
$\nest\Phi$ consists only of the empty word if $\Phi$ is atomic.
\item
$\nest{\neg\Phi}=\overline{\nest\Phi}$.
\item
$\nest{\Phi\And\Psi}=\nest{\Phi\Or\Psi}=\nest\Phi\cup\nest\Psi$.
\item
$\nest{\exists x\Phi}=\exists\nest\Phi$ and
$\nest{\forall x\Phi}=\forall\nest\Phi$.
\end{enumerate}
\end{definition}

\begin{definition}
The {\em quantifier rank\/} of a formula $\Phi$
is the maximum length of a string in $\nest\Phi$.
\end{definition}

\noindent
We adopt the notion of the {\em alternation number\/} of a formula
(cf.\ \cite[Definition 2.8]{Pez1}).

\begin{definition}
Given a sequence of quantifiers $s$, let $\alt s$ denote the number
of occurrences of $\exists\forall$ and $\forall\exists$ in $s$.
The {\em alternation number\/} of a first order formula $\Phi$
is the maximum $\alt s$ over $s\in\nest\Phi$.
\end{definition}

\begin{definition}
Given a graph $G$ and a first order formula $\Phi$ over vocabulary
$\{E,{=}\}$, we say that $\Phi$ {\em distinguishes $G$ from $G'$} if
$\Phi$ is true on $G$ but false on $G'$.
By $\D{G,G'}$ (resp.\ $\DD k{G,G'}$) we denote the minimum quantifier rank of
a formula (with alternation number at most $k$ resp.) distinguishing $G$
from $G'$.

We say that $\Phi$ {\em defines $G$\/} (up to isomorphism) if $\Phi$
distinguishes $G$ from any non-isomorphic graph $G'$.
By $\D G$ (resp.\ $\DD kG$) we denote the minimum
quantifier rank of a formula defining $G$
(with alternation number at most $k$ resp.).
\end{definition}

Note that $\D{G,G'}=\min_{k\ge0}\DD k{G,G'}$ and $\D G=\min_{k\ge0}\DD kG$.
The following proposition is a simple consequence of the
well-known fact that over a fixed finite relational vocabulary
there are only finitely many inequivalent first order formulas of
bounded quantifier rank.

\begin{proposition}\label{prop:ddd}
For every graph $G$ it holds
$\D G = \max\setdef{\D{G,G'}}{G'\not\cong G}$ and
$\DD kG = \max\setdef{\DD k{G,G'}}{G'\not\cong G}$.
\end{proposition}

\subsection{Games}

The {\em \EF\/ game\/} is played on a pair of structures of the same
type. We give a definition in the terminology of graphs.

\begin{definition}
Let $G$ and $G'$ be graphs with disjoint vertex sets.
The $r$-round \EF\/ game on $G$ and $G'$,
denoted by $\game_r(G,G')$, is played by
two players, Spoiler and Duplicator, with using $r$ pairwise distinct
pebbles $p_1,\ldots,p_r$, each given in duplicate. Spoiler starts the game.
A {\em round\/} consists of a move of Spoiler followed by a move of
Duplicator. In the $s$-th round Spoiler selects one of
the graphs $G$ or $G'$ and places $p_s$ on a vertex of this graph.
In response Duplicator should place the other copy of $p_s$ on a vertex
of the other graph. It is allowed to place more than one pebble on the
same vertex. We will use $x_s$ (resp.\ $y_s$) to
denote the vertex of $G$ (resp.\ $G'$) occupied by $p_s$, irrespectively
of who of the players places the pebble on this vertex.
The pair of sequences $\bar x=(x_1,\ldots,x_s)$ and
$\bar y=(y_1,\ldots,y_s)$ is a {\em configuration\/} of the game
after the $s$-th round.
If after every of $r$ rounds it is true that
$$
x_i=x_j\mbox{ iff } y_i=y_j\mbox{ for all }i,j\le s,
$$
and the component-wise correspondence $\bar x$ to
$\bar y$ is a partial isomorphism from $G$ to $G'$, this is
a win for Duplicator;  Otherwise the winner is Spoiler.

The {\em $k$-alternation\/} \EF\/ game on $G$ and $G'$ is a variant
of the game in which Spoiler is allowed to switch from one graph
to another at most $k$ times during the game, i.e., in at most $k$
rounds he can choose the graph other than that in the preceding round.
\end{definition}

\noindent
The main technical tool we will use is given by the following statement.

\begin{proposition}\label{prop:loggames}
Let $G$ and $G'$ be non-isomorphic graphs.

\begin{enumerate}
\item
$\D{G,G'}$ equals the minimum $r$ such that Spoiler has a winning
strategy in $\game_r(G,G')$.
\item
$\DD k{G,G'}$ equals the minimum $r$ such that Spoiler has a winning
strategy in the $k$-alternation $\game_r(G,G')$.
\end{enumerate}
\end{proposition}

\noindent
We refer the reader to \cite[Theorem 1.2.8]{EFl}, \cite[Theorem 6.10]{Imm},
or \cite[Theorem 2.3.1]{Spe} for the proof of the first
claim and to \cite{Pez2} for the second claim.

\subsection{Notation}

Throughout the paper $\log n$ means the logarithm base 2.
We will use $[n]$ to denote the set $\{1,\ldots,n\}$.
For a function $f$, by $f^{(k)}$ we will denote the $k$-fold composition
of~$f$.

\section{Spoiler's winning strategy based on separators}\label{s:str}

\begin{lemma}\label{lem:conn}
Let $\bar x$, $\bar y$ be a configuration after the $s$-th round
of the Ehrenfeucht game on colored graphs $G$, $G'$.
Given $I\subset[s]$, let $X_I=\setdef{x_h}{h\in I}$ and
$Y_I=\setdef{y_h}{h\in I}$.
Assume that there are $i,j\in[s]$ and $I\subset[s]$ such that
$x_i$ and $x_j$ are in the same $X_I$-flap $F$ of $G$ while
$y_i$ and $y_j$ are in different $Y_I$-flaps of $G'$.
Then Spoiler is able to win in at most $\lceil\log|F|\rceil$ moves,
playing only in $G$.
\end{lemma}

\begin{proof}
Let $d_F$ be the standard metric on $F$, i.e., $d_F(z_1,z_2)$ is
the shortest length of a path from $z_1$ to $z_2$ through vertices of~$F$.

Spoiler sets $u_1=x_i$, $u_2=x_j$, $v_1=y_i$, $v_2=y_j$, and places
a pebble on a vertex $u$ in $F$ such that
$d_F(u,u_m)\le\lceil d_F(u_1,u_2)/2\rceil$ for both $m=1,2$.
Let $v\in V(G')$ be selected by Duplicator in response to $u$.
For both $m=1,2$ there is a path of length at most
$\lceil d_F(u_1,u_2)/2\rceil$ from $u$ to $u_m$ through vertices
of $G\setminus X_I$. In contrast with this, for $m=1$ or $m=2$
there is no path from $v$ to $v_m$ through vertices of
$G'\setminus Y_I$. For this particular value of $m$,
Spoiler resets $u_1=u$, $u_2=u_m$, $v_1=v$, $v_2=v_m$ and applies
the same strategy once again. Therewith Spoiler ensures that, in each round,
$v_1$ and $v_2$ are in different $Y_I$-flaps of $G'$. Eventually,
unless Duplicator loses earlier, $d_F(u,u_m)=1$ for both $m=1,2$, that is,
$u_1,u,u_2$ is a path in $F$. As there is no path from $v_1$ to $v_2$
inside $G'\setminus Y_I$, the isomorphism is broken and Spoiler wins.

To estimate the number of moves made, notice that initially
$d_F(u_1,u_2)\le|F|-1$ and for each subsequent $u_1,u_2$ this distance
becomes at most
$f(d_F(u_1,u_2))$, where $f(\ga)=(\ga+1)/2$. Therefore the number of
moves does not exceed the minimum $k$ such that $f^{(k)}(|F|-1)<2$.
As $(f^{(k)})^{-1}(\gb)=2^k\gb-2^k+1$, the latter inequality is equivalent to
$2^k>|F|-2$, which proves the bound.
\end{proof}

\begin{remark}
The bound $\lceil\log|F|\rceil$ in Lemma \ref{lem:conn} is tight up to
a small additive constant. For example, let $2C_n$ be the disjoint union
of two cycles of length $n$. It is known (e.g.\ \cite[Example 2.3.8]{EFl})
that Duplicator can survive in the \EF\/ game on $2C_n$ and $C_n$
during $\lfloor\log(n-1)\rfloor$ rounds for any strategy of Spoiler,
in particular, when Spoiler's first move is in one component of $2C_n$
and his second move is in the other component of $2C_n$
($C_n$ and both components of $2C_n$ are considered $\emptyset$-flaps).
\end{remark}

We now give a fairly rudimentary definition of a separator that
abstracts away some important features usually associated
with this notion. These features will be specified a bit later
(cf.\ Definition \ref{def:sep}).

\begin{definition}
Let $\cl$ be a class of graphs. Let $k$ be a function from $\bbn$ to $\bbr$.
A {\em separator of size $k$\/} for $\cl$ is a function $\Sep$ defined on
graphs in $\cl$ such that, for every $G\in\cl$, $\sep G\subseteq V(G)$
and $|\sep G|\le k(n)$. We call the separator {\em hereditary\/} if,
for every $G\in\cl$, all $\sep G$-flaps are in~$\cl$.
We will extend $\Sep$ over colored graphs by setting
$\sep{\tilde G}=\sep G$ for an arbitrary coloration $\tilde G$
of a graph~$G$.
\end{definition}

Let $\cl$ be a class of connected graphs with hereditary separator $\Sep$ of
size $k(n)$. For each $i\ge 0$, we now describe a winning strategy $S_i$
for Spoiler applicable to the \EF\/ game on colored graphs $G$ and $G'$
such that $G\in\col\cl$ and $G'$ is arbitrary non-isomorphic with $G$.
The strategy $S_i$ is designed so that Spoiler reduces play on $G$
to play on a $\sep G$-flap $F$, where the index $i$ means that
this trick can be used up to $i$ times.

\medskip

{\sc Strategy $S_0$.}

\smallskip

\noindent
Spoiler selects all vertices of $G$. If this is still not a win
($G$ is isomorphic to a proper induced subgraph of $G'$), then
Spoiler selects one more vertex $u$ in $G'$. This is Spoiler's win
for any choice of $u$. Nevertheless, we will need an additional
condition: If possible (e.g., $G'$ is also connected),
$u$ must be adjacent to a vertex that was selected
in $G'$ by Duplicator.

\medskip

{\sc Strategy $S_i$, $i\ge 1$.}

\smallskip

\noindent
As usually, we denote the order of $G$ by $n$.
If $G'$ is disconnected, Spoiler selects 2 vertices in different
components of $G'$ and then applies the strategy of Lemma \ref{lem:conn}
(with $I=\emptyset$) winning in at most $\lceil\log n\rceil$ next moves.
Assume that $G'$ is connected.

If $k(n)\ge n$, then Spoiler applies the strategy $S_0$.
Assume that $k(n)<n$. Let $X=\sep G$ and $k=|X|$. Spoiler first selects
all $k$ vertices of $X$. Let $(x_1,\ldots,x_k)$, $(y_1,\ldots,y_k)$
be the configuration at this stage of the game. Note that
$X=\{x_1,\ldots,x_k\}$ and let $X'=\{y_1,\ldots,y_k\}$.
Assume that Duplicator still does not lose.

We now modify coloring of each $X$-flap and each $X'$-flap in the following
way. Let $A_1,\ldots,A_k$ be colors that occur neither in $G$ nor in $G'$.
If $F$ is $X$-flap (resp.\ $X'$-flap), then, for every $i\le k$, the color
$A_i$ is assigned to all those vertices of $F$ that are adjacent to $x_i$
in $G$ (resp.\ to $y_i$ in $G'$).
It should be stressed that we introduce new colors solely for the sake
of technical convenience and this puts no new constraint on Duplicator:
If Duplicator violates a new color, he therewith violates the adjacency
to a previously selected vertex. Our technical benefit is that now,
if Spoiler and Duplicator play on an $X$-flap $F$ and $X'$-flap $F'$,
we can forget about the rest of $G$ and $G'$ because Spoiler's win
in the game on $F$ and $F'$ will mean his win in the game on $G$ and $G'$
as well.

Given a colored graph $H$, let $m(H)$ (resp.\ $m'(H)$) denote
the number of $X$-flaps (resp.\ $X'$-flaps) isomorphic to $H$
(the flaps are assumed with modified colorings).
Observe that the partial isomorphism $\phi\function X{X'}$
established by Duplicator extends to an isomorphism from $G$ to $G'$
iff $m(H)=m'(H)$ for all $H$. As $G\not\cong G'$, there is $H$
with $m(H)\ne m'(H)$.

{\it Case 1:
There is an $H$ with $m(H)>m'(H)$.}
Spoiler starts to select one vertex in each $X$-flap isomorphic to $H$.
At latest in the $(m'(H)+1)$-th round Duplicator either selects
two vertices in the same $X'$-flap isomorphic to $H$ or selects a vertex
in an $X'$-flap $F'$ non-isomorphic to $H$. In the former case Spoiler
applies the strategy of Lemma \ref{lem:conn} and wins in at most
$\lceil\log|H|\rceil$ moves. In the latter case Spoiler applies
the strategy $S_{i-1}$ on graphs $F$ and $F'$, where $F$ is the
$X$-flap, isomorphic to $H$, that is visited by Spoiler in the same
round when Duplicator visits $F'$. Spoiler follows $S_{i-1}$ on $F$ and $F'$
as long as Duplicator plays inside $F$ and $F'$. If in any round
Spoiler selects a vertex in $F$ but Duplicator replies with a vertex
outside $F'$, Spoiler switches to the strategy of Lemma \ref{lem:conn}
and wins in at most $\lceil\log|F|\rceil$ moves. We do not need to specify
the strategy $S_i$ for the case when Spoiler selects a vertex in $F'$
but Duplicator replies with a vertex outside $F$ because this is a loss
of Duplicator (see Lemma \ref{lem:out} below).
Note that we so far have not encountered the situation when Spoiler
switches from $G$ to $G'$, except the strategy~$S_0$.

{\it Case 2:
$m(H)\le m'(H)$ for all $H$.}
It follows that the number of $X$-flaps, that will be denoted by $f$,
is strictly less than the number of $X'$-flaps. In the subsequent
$f$ moves Spoiler selects one vertex in each $X$-flap.
Then, in the next move Spoiler selects a vertex $u$ in an $X'$-flap
$F'$ that was not visited by Duplicator so that $u$ is adjacent to
a vertex in $X'$ (such a vertex exists in $F'$ because $G'$ is connected
and there can be no path from $F'$ to $X'$ through other flaps).
Duplicator is enforced to select in response yet another, the second,
vertex in one of $X$-flaps, say, in $F$. Then Spoiler applies the
strategy of Lemma \ref{lem:conn} and wins in at most $\lceil\log|F|\rceil$
moves. This completes the description of the strategy $S_i$.

\begin{lemma}\label{lem:out}
Let $\cl$ be a class of connected graphs with hereditary separator $\Sep$.
Let $t\ge 1$. Consider the \EF\/ game on colored graphs $G\in\col\cl$ and
$G'\not\cong G$ in
which Spoiler follows $S_t$ and Duplicator follows an arbitrary fixed
strategy. Assume that during the course of the game Spoiler
invokes the strategy $S_{t-1}$ on $X$-flap $F$ and $X'$-flap $F'$.
If afterwards Spoiler selects a vertex in $F'$ but Duplicator
responds with a vertex outside $F$, this is an immediate loss for
Duplicator.
\end{lemma}

\begin{proof}
We start with some notation.
Let $q$, $0\le q\le t-1$, be the smallest number such that
the strategy $S_q$ is invoked during the course of the game.
Set $F_t=G$, $F'_t=G'$, and $X_{t+1}=X'_{t+1}=\emptyset$.
If $q\le i\le t$, let $X_i=\sep{F_i}$ and $X'_i$ be the subset of
$V(G')$ that Duplicator takes into correspondence with $X_i$.
Furthermore, let $F_i$ and $F'_i$ be respectively $X_{i+1}$-
and $X'_{i+1}$-flaps on which Spoiler applies the strategy $S_i$.
In particular, $X=X_t$, $X'=X'_t$, $F=F_{t-1}$, and $F'=F'_{t-1}$.
Note that $X_i\subset V(F_i)$ and $V(F_i)\subset V(F_{i+1})$.

Tracing trough the description of $S_t$, it is easy to see that,
as $S_{t-1}$ is known to be invoked, Spoiler can change $G$ for $G'$
only when invoking $S_l$ for some $q\le l\le t-1$ and, moreover,
if $l>0$, there must be Case 2.

Assume that Spoiler selects a vertex $u$ in $F'$ and consider
first the case that $l>0$. Thus, $u$ is selected in the $X'_l$-flap
$F'_{l-1}$ so that $d(u,X'_l)=1$. Assume that Duplicator responds
with $v\in V(G)$ and does not lose. Therefore $d(v,X_l)=1$ too.
Notice that $\Gg(X_l)\subseteq V(F_l)\cup\bigcup^t_{i=l+1}X_i$.
Since all vertices in the $X_i$'s are already occupied, we conclude
that $v\in V(F_l)\subset\ldots\subset V(F_{t-1})=V(F)$.

Consider now the case that $l=0$. Spoiler selects $u$ in the last round
of the play on $F_0$ and $F'_0$ with the strategy $S_0$.
Recall that $u$ is selected
in $F'_0$ so that $u$ is adjacent to a vertex $y$ previously selected
in $F'_0$ be Duplicator. Let $x$ be a counterpart of $y$ in $F_0$.
Hence Duplicator is forced to select a $v\in V(G)$ that is adjacent to
$x$. Since $\Gg(F_0)\subseteq\bigcup^t_{i=1}X_i$ consists of previously
occupied  vertices, Duplicator loses anyway.
\end{proof}

\begin{definition}\label{def:sep}
Let $m$ be an integer and $\eps\in(0,1)$. A separator $\Sep$ for
a class of graphs $\cl$ is called an {\em $m$-flap $\eps$-separator\/}
if, for every $G\in\cl$, there are at most $m$ $\sep G$-flaps
and each of them has at most $\eps n$ vertices.
\end{definition}

We now estimate the length of the \EF\/ game if Spoiler follows
the strategy $S_t$, where the choice of $t$ is optimized for our purposes.

\begin{lemma}\label{lem:length}
Suppose that $\cl$ is a class of connected graphs with hereditary $m$-flap
$\eps$-separator $\Sep$ of size $k(n)$, where the function $k$ is
assumed to be defined over reals and non-decreasing.
For $i\ge 0$, let the strategies $S_i$ be based on $\Sep$.
If $G$ and $G'$ are colored graphs such that $G\in\col\cl$ and
$G'\not\cong G$, let $L_i(G,G')$ denote the minimum $r$ such that
Spoiler wins $\game_r(G,G')$ following $S_i$, regardless of Duplicator's
strategy. Set $t=\lceil\log(n/m)/\log(\eps^{-1})\rceil$. Then
$$
L_t(G,G')<\sum_{i=0}^{t-1}k(\eps^i n)+m(t+1)+\log n+2.
$$
\end{lemma}

\begin{proof}
Assume that $H\in\col\cl$ and $H'$ is an arbitrary connected colored
graph. Fix the strategy $S_i$ for Spoiler and an arbitrary strategy
$D$ for Duplicator. Consider the \EF\/ game on $H$ and $H'$ in which
the players follow these strategies. Let $l_i(H,H')$ denote the number
of rounds till Spoiler wins or a position as in Lemma \ref{lem:conn}
occurs. Note first that
\begin{equation}\label{eq:hh0}
l_0(H,H')\le |H|+1.
\end{equation}
Let $n=|H|$. If $n$ is such that $k(n)\ge n$, we have
\begin{equation}\label{eq:hh1}
l_i(H,H')=l_0(H,H')\le n+1\le k(n)+1.
\end{equation}
Otherwise Spoiler selects all vertices of $\sep H$ and the further play
depends on which of Cases 1 or 2 takes place. In Case 2 we have
\begin{equation}\label{eq:hh2}
l_i(H,H')\le k(n)+m+1.
\end{equation}
In Case 1, suppose that the strategy $S_{i-1}$ is invoked on subgraphs
$F$ and $F'$. Then it is easy to see that
\begin{equation}\label{eq:hh3}
l_i(H,H')\le k(n)+m+l_{i-1}(F,F')
\end{equation}
and this bound exceeds the bounds \refeq{eq:hh1} and \refeq{eq:hh2}.
Notice that $|F|\le\eps n$.

Assume that $G'$ is connected and let $n=|G|$.
By a simple inductive argument, \refeq{eq:hh3} and \refeq{eq:hh0}
imply that
$$
l_t(G,G')\le \sum_{i=0}^{t-1}k(\eps^i n)+mt+\eps^t n+1\le
\sum_{i=0}^{t-1}k(\eps^i n)+m(t+1)+1,
$$
where the latter inequality follows by the choice of $t$.
Since Duplicator's strategy $D$ was chosen arbitrary, by the description
of $S_i$ and Lemma \ref{lem:conn} we have
$$
L_t(G,G')\le l_t(G,G')+\lceil\log(\eps n)\rceil<
\sum_{i=0}^{t-1}k(\eps^i n)+m(t+1)+\log n+2.
$$

In the case that $G'$ is disconnected, we easily have a better bound
$$
L_t(G,G')\le \lceil\log n\rceil + 2.
$$
This completes the proof.
\end{proof}

\begin{lemma}\label{lem:const}
Suppose that $\cl$ is a class of connected graphs with hereditary $m$-flap
$\eps$-separator of size $k$, where $k$ is a constant and $\eps m > 1$.
Then for every $G\in\cl$ it holds
$$
\DD 2 G < \of{\frac{k+m}{\log(\eps^{-1})}+1}\log n + m + 2.
$$
\end{lemma}

\begin{proof}
By Proposition \ref{prop:ddd}, we actually have to estimate
$\DD2{G,G'}$ uniformly for all $G'\not\cong G$. We use the
characterization of $\DD2{G,G'}$ as the length of the \EF\/ game
given by Proposition \ref{prop:loggames}. We may consider $G$ and $G'$
colored graphs (whose vertices satisfy no color relation), because
the length of the game will obviously remain the same. As $G\in\cl$,
Spoiler can apply the strategy $S_t$ with $t$ as in Lemma \ref{lem:length}.
By this lemma, he wins in less than $tk+m(t+1)+\log n+2$ moves with
$$
t<\frac{\log(n/m)}{\log(\eps^{-1})}+1<\frac{\log n}{\log(\eps^{-1})},
$$
the latter inequality because $\eps m > 1$.
This easily implies the bound claimed.
It remains to notice that, following $S_t$, Spoiler makes at most 2
alternations between $G$ and $G'$. The alternation from $G$ to $G'$
happens when, under the invocation of some $S_l$, $l\le t$, for the
first time there occurs Case 2. Afterwards Spoiler goes back to
$G$ and wins without further invocation of $S_{l-1}$ and hence
without further alternation.
\end{proof}

\section{Application to trees and outerplanar graphs}\label{s:hop}

\begin{theorem}\label{thm:trees}
Let $d\ge 2$ be a constant. If $G$ is a tree of maximum vertex degree $d$,
then
\begin{equation}\label{eq:trees}
\DD2G < c_d\log n + d + 2,
\end{equation}
where $c_d=(d+1)/\log(3/2)+1$.
\end{theorem}

\begin{proof}
We use the well-known fact that every tree
has $\frac23$-separator of size 1 (see e.g.\ \cite{LTa}
Let $\cl_d$ consists of all trees of degree at most $d$.
Any separator for $\cl_d$ is obviously hereditary, and it is easy to
see that a 1-vertex separator splits a $G\in\cl_d$ in at most $d$
flaps. Lemma \ref{lem:const} applies to $\cl_d$ and implies~\refeq{eq:trees}.
\end{proof}

\begin{remark}\label{rem:path}
The constants $c_d$ in Theorem \ref{thm:trees} are not the best possible.
For example, if $d=2$, the optimal constant is 1 because
$\DD1{P_n}<\log n+3$. Indeed, if $G$ is disconnected, then
$\DD1{P_n,G}<\log n+3$ by Lemma \ref{lem:conn}. If $G$ is connected
with $\Delta(G)>2$, then clearly $\DD0{P_n,G}\le4$. If $G=C_m$, then
$\DD1{P_n,G}\le3$. Finally, if $G=P_m$ is another path, then
$\DD1{P_n,G}<\log n+3$ by \cite[Theorem 2.1.2]{Spe}.
\end{remark}

\begin{theorem}\label{thm:hop}
If $G$ is a Hamiltonian outerplanar graph, then
$$
\DD2G < c\log n + 9,
$$
where $c=12/\log(3/2)+1<22$.
\end{theorem}

The theorem directly follows from Lemma \ref{lem:const}
and Lemma \ref{lem:hop} below. It is well known \cite{Lei} that
every outerplanar graph has $\frac23$-separator of size 2.
If a graph is in addition Hamiltonian, then there are obviously only
2 flaps. However, for the class of Hamiltonian outerplanar graphs
such separator is not hereditary. Fortunately, we are able to
extend the class of Hamiltonian outerplanar graphs to a wider class
with hereditary constant-flap $\frac23$-separator of constant size.

\begin{definition}
We will further on abbreviate the term {\em Hamiltonian outerplanar\/}
as \hop. An {\em edge-deleted \hop\/} graph (to be abbreviated as \edhop)
is a graph which is not \hop\/ but becomes such after joining a pair of
vertices by an edge. A {\em 2-edge-deleted \hop\/} graph (to be abbreviated
as \eedhop) is a connected graph which is not \edhop\/ but becomes
such after joining a pair of vertices by an edge.
Let $\clo$ be the class of \hop, \edhop, and \eedhop\/ graphs.
\end{definition}

\begin{lemma}\label{lem:hop}
The class $\clo$ has hereditary 7-flap $\frac23$-separator
of size~5.
\end{lemma}

\begin{proof}
We will use the fact that outerplanar graphs have 2-vertex
$\frac23$-separators. Assume that $G\in\clo$ and show the existence
of an appropriate separator in~$G$. Provided any separator $X$ of $G$
is specified, by $\flap{v_1,\ldots,v_q}$ we will denote the $X$-flap
containing the vertices $v_1,\ldots,v_q$ (if such exists).

{\it Case 1: $G$ is \hop.}
Let $\{s_1,s_2\}$ be a $\frac23$-separator of $G$ and $C$ be a Hamilton
cycle in $G$. For $i=1,2$, denote the neighbors of $s_i$ in $C$ by
$a_i$ and $b_i$ as shown in Figure 1. The separator splits $G$
in two flaps $\flap{a_1,a_2}$ and $\flap{b_1,b_2}$.
If $a_1$ and $a_2$ (resp.\ $b_1$ and $b_2$) are adjacent, the flap
$\flap{a_1,a_2}$ (resp.\ $\flap{b_1,b_2}$)
is \hop, otherwise it is \edhop.
For the completeness notice the possibility that $a_1=a_2$ or $b_1=b_2$,
which can happen for graphs on at most 6 vertices.

\begin{figure}
\centerline{
\unitlength=1.00mm
\special{em:linewidth 0.4pt}
\linethickness{0.4pt}
\begin{picture}(143.00,95.00)
\put(22.00,81.00){\oval(14.00,20.00)[l]}
\put(37.00,81.00){\oval(14.00,20.00)[r]}
\put(22.00,71.00){\circle*{2.00}}
\put(22.00,91.00){\circle*{2.00}}
\put(37.00,91.00){\circle*{2.00}}
\put(37.00,71.00){\circle*{2.00}}
\put(30.00,92.00){\circle*{2.00}}
\put(30.00,70.00){\circle*{2.00}}
\put(22.00,94.00){\makebox(0,0)[cb]{$a_1$}}
\put(37.00,94.00){\makebox(0,0)[cb]{$b_1$}}
\put(30.00,95.00){\makebox(0,0)[cb]{$s_1$}}
\put(22.00,68.00){\makebox(0,0)[ct]{$a_2$}}
\put(37.00,68.00){\makebox(0,0)[ct]{$b_2$}}
\put(30.00,67.00){\makebox(0,0)[ct]{$s_2$}}
\put(99.00,81.00){\oval(14.00,20.00)[l]}
\put(99.00,71.00){\circle*{2.00}}
\put(99.00,91.00){\circle*{2.00}}
\put(114.00,91.00){\circle*{2.00}}
\put(114.00,71.00){\circle*{2.00}}
\put(107.00,92.00){\circle*{2.00}}
\put(107.00,70.00){\circle*{2.00}}
\put(99.00,94.00){\makebox(0,0)[cb]{$a_1$}}
\put(114.00,94.00){\makebox(0,0)[cb]{$b_1$}}
\put(107.00,95.00){\makebox(0,0)[cb]{$s_1$}}
\put(99.00,68.00){\makebox(0,0)[ct]{$a_2$}}
\put(114.00,68.00){\makebox(0,0)[ct]{$b_2$}}
\put(107.00,67.00){\makebox(0,0)[ct]{$s_2$}}
\put(114.00,74.50){\oval(12.00,7.00)[rb]}
\put(113.00,86.00){\oval(14.00,10.00)[rt]}
\put(120.00,86.00){\circle*{2.00}}
\put(120.00,75.00){\circle*{2.00}}
\put(120.00,75.00){\dashbox{1.00}(0.00,11.00)[cc]{}}
\put(123.00,86.00){\makebox(0,0)[lc]{$c_1$}}
\put(123.00,75.00){\makebox(0,0)[lc]{$c_2$}}
\emline{10.00}{35.00}{1}{20.00}{37.00}{2}
\emline{40.00}{37.00}{3}{50.00}{35.00}{4}
\emline{10.00}{20.00}{5}{20.00}{18.00}{6}
\emline{40.00}{18.00}{7}{50.00}{20.00}{8}
\put(10.00,35.00){\circle*{2.00}}
\put(20.00,37.00){\circle*{2.00}}
\put(30.00,38.00){\circle*{2.00}}
\put(40.00,37.00){\circle*{2.00}}
\put(50.00,35.00){\circle*{2.00}}
\put(50.00,20.00){\circle*{2.00}}
\put(40.00,18.00){\circle*{2.00}}
\put(30.00,17.00){\circle*{2.00}}
\put(20.00,18.00){\circle*{2.00}}
\put(10.00,20.00){\circle*{2.00}}
\put(10.00,20.00){\dashbox{1.00}(0.00,15.00)[cc]{}}
\put(50.00,20.00){\dashbox{1.00}(0.00,15.00)[cc]{}}
\put(7.00,35.00){\makebox(0,0)[rc]{$d_1$}}
\put(7.00,20.00){\makebox(0,0)[rc]{$d_2$}}
\put(53.00,35.00){\makebox(0,0)[lc]{$c_1$}}
\put(53.00,20.00){\makebox(0,0)[lc]{$c_2$}}
\put(20.00,40.00){\makebox(0,0)[cb]{$a_1$}}
\put(30.00,41.00){\makebox(0,0)[cb]{$s_1$}}
\put(40.00,40.00){\makebox(0,0)[cb]{$b_1$}}
\put(20.00,15.00){\makebox(0,0)[ct]{$a_2$}}
\put(30.00,14.00){\makebox(0,0)[ct]{$s_2$}}
\put(40.00,15.00){\makebox(0,0)[ct]{$b_2$}}
\emline{73.00}{35.00}{9}{88.00}{41.00}{10}
\emline{100.00}{43.00}{11}{113.00}{43.00}{12}
\emline{125.00}{41.00}{13}{140.00}{35.00}{14}
\put(73.00,35.00){\circle*{2.00}}
\put(88.00,41.00){\circle*{2.00}}
\put(100.00,43.00){\circle*{2.00}}
\put(113.00,43.00){\circle*{2.00}}
\put(125.00,41.00){\circle*{2.00}}
\put(140.00,35.00){\circle*{2.00}}
\put(78.00,37.00){\circle*{2.00}}
\put(83.00,39.00){\circle*{2.00}}
\emline{83.00}{39.00}{15}{134.00}{37.00}{16}
\put(134.00,37.00){\circle*{2.00}}
\put(94.00,42.00){\circle*{2.00}}
\put(119.00,42.00){\circle*{2.00}}
\emline{73.00}{20.00}{17}{107.00}{17.00}{18}
\emline{107.00}{17.00}{19}{140.00}{20.00}{20}
\put(73.00,20.00){\circle*{2.00}}
\put(107.00,17.00){\circle*{2.00}}
\put(140.00,20.00){\circle*{2.00}}
\emline{78.00}{37.00}{21}{107.00}{17.00}{22}
\put(73.00,20.00){\dashbox{1.00}(0.00,15.00)[cc]{}}
\put(140.00,20.00){\dashbox{1.00}(0.00,15.00)[cc]{}}
\put(70.00,35.00){\makebox(0,0)[rc]{$d_1$}}
\put(76.00,39.00){\makebox(0,0)[cb]{$e_1$}}
\put(100.00,46.00){\makebox(0,0)[cb]{$b_1$}}
\put(113.00,46.00){\makebox(0,0)[cb]{$a_2$}}
\put(134.00,40.00){\makebox(0,0)[lb]{$f$}}
\put(143.00,35.00){\makebox(0,0)[lc]{$c_1$}}
\put(143.00,20.00){\makebox(0,0)[lc]{$c_2$}}
\put(70.00,20.00){\makebox(0,0)[rc]{$d_2$}}
\put(107.00,14.00){\makebox(0,0)[ct]{$e_2$}}
\put(30.00,60.00){\makebox(0,0)[ct]{Case 1}}
\put(107.00,60.00){\makebox(0,0)[ct]{Case 2}}
\put(30.00,7.00){\makebox(0,0)[ct]{Subcase 3.1}}
\put(107.00,7.00){\makebox(0,0)[ct]{Subcase 3.2}}
\put(88.00,44.00){\makebox(0,0)[cb]{$a_1$}}
\put(94.00,45.00){\makebox(0,0)[cb]{$s_1$}}
\put(119.00,45.00){\makebox(0,0)[cb]{$s_2$}}
\put(125.00,44.00){\makebox(0,0)[cb]{$b_2$}}
\bezier{20}(22.00,91.00)(30.00,93.00)(37.00,91.00)
\bezier{20}(22.00,71.00)(30.00,69.00)(37.00,71.00)
\bezier{20}(99.00,91.00)(107.00,93.00)(114.00,91.00)
\bezier{20}(99.00,71.00)(107.00,69.00)(114.00,71.00)
\bezier{30}(20.00,37.00)(30.00,39.00)(40.00,37.00)
\bezier{30}(20.00,18.00)(30.00,16.00)(40.00,18.00)
\bezier{20}(88.00,41.00)(94.00,42.00)(100.00,43.00)
\bezier{20}(113.00,43.00)(119.00,42.00)(125.00,41.00)
\end{picture}
}
\caption{Proof of Lemma \protect\ref{lem:hop}}
\end{figure}

{\it Case 2: $G$ is \edhop.}
Let $\{c_1,c_2\}$ be the edge that $G$ lacks to become a \hop\/ graph
$\bar G$. Let $C$ be a Hamilton
cycle of $\bar G$. Given $u,v\in V(G)$, by $[u,v]$ we will denote
the set of vertices on the arc of $C$ from $u$ to $v$ that does not contain
$\{c_1,c_2\}$ (if such exists).

Let $\{s_1,s_2\}$ be a $\frac23$-separator of $\bar G$.
If $\{s_1,s_2\}$ and $\{c_1,c_2\}$ intersect,
$\{s_1,s_2\}$-flaps of $G$ and $\bar G$ are the same, and we essentially
have Case 1. Assume that $\{s_1,s_2\}$ and $\{c_1,c_2\}$ are disjoint
and denote the neighbors of $s_i$ in $C$ by $a_i$ and $b_i$ as in Figure 1
(some of the vertices shown there may coincide).

Note that
$[a_1,a_2]$ can be connected neither to $[b_1,c_1]$ nor to $[b_2,c_2]$.
As easily seen, the flap $\flap{a_1,a_2}$ is either \hop\/ or \edhop,
depending on if $a_1$ and $a_2$ are adjacent.
If $[b_1,c_1]$ and $[b_2,c_2]$ are connected by an edge,
$\flap{b_1,c_1}=\flap{b_2,c_2}$ is the second $\{s_1,s_2\}$-flap
which is \edhop\/ or \eedhop, depending on the adjacency of $b_1$ and $b_2$;
otherwise $\flap{b_1,c_1}$ and $\flap{b_2,c_2}$ are
two flaps, \hop\/ or \edhop\/ depending on the adjacency of
$b_1$, $c_1$ and $b_2$, $c_2$ respectively.

{\it Case 3: $G$ is \eedhop.}
Let $\{c_1,c_2\}$ and $\{d_1,d_2\}$ be the edges that $G$ lacks to be
a \hop\/ graph $\bar G$. Let $C$ be a Hamilton
cycle of $\bar G$. Given $u,v\in V(G)$, by $[u,v]$ we now denote
the set of vertices on the arc of $C$ from $u$ to $v$ that does not contain
$\{c_1,c_2\}$ and $\{d_1,d_2\}$ (if such exists). Furthermore,
$[u,v)=[u,v]\setminus\{v\}$ and $(u,v]=[u,v]\setminus\{u\}$.

Let $\{s_1,s_2\}$ be a $\frac23$-separator of $\bar G$.
If $\{s_1,s_2\}$ and $\{c_1,c_2,d_1,d_2\}$ intersect, we essentially
have Case 2. Assume they do not.
Denote the neighbors of $s_i$ in $C$ by $a_i$ and $b_i$.
Two different subcases, depicted in Figure 1, are possible
(some of the vertices may coincide).

{\it Subcase 3.1: $s_1$ and $s_2$ are separated in $C$ by the missing
edges $\{c_1,c_2\}$ and $\{d_1,d_2\}$.}
Note that neither $[d_1,a_1]$ nor $[d_2,a_2]$ can be connected to
any of $[b_1,c_1]$ and $[b_2,c_2]$. If $[d_1,a_1]$ is connected to
$[d_2,a_2]$, the flap $\flap{d_1,a_1}=\flap{d_2,a_2}$ is \edhop\/
or \eedhop, depending on whether $a_1$ and $a_2$ are adjacent or not.
Otherwise, $\flap{d_1,a_1}$ and $\flap{d_2,a_2}$ are two flaps,
each \hop\/ or \edhop. All the same is true for the pair
$[b_1,c_1]$ and $[b_2,c_2]$.

{\it Subcase 3.2: $s_1$ and $s_2$ are in $C$ not separated by
$\{c_1,c_2\}$ and $\{d_1,d_2\}$.}
Note that $[b_1,a_2]$ cannot be connected to any of
$[d_1,a_1]$, $[b_2,c_1]$, and $[d_2,c_2]$, and that
$\flap{b_1,a_2}$ is a \hop\/
or \edhop\/ flap. If one of $[d_1,a_1]$, $[b_2,c_1]$, and $[d_2,c_2]$
is disconnected to the other two, then the corresponding of the components
$\flap{d_1,a_1}$, $\flap{b_2,c_1}$, and $\flap{d_2,c_2}$ is \hop\/
or \edhop, and the other two are also \hop\/ or \edhop\/ (if there is no
edge between the corresponding arcs) or they are the same \edhop\/ or
\eedhop\/ flap (if a such edge exists).

Consider now the case that one of $[d_1,a_1]$, $[b_2,c_1]$, and $[d_2,c_2]$
is connected to the other two, that is,
$\flap{d_1,a_1}=\flap{b_2,c_1}=\flap{d_2,c_2}$.
Without loss of generality, assume that $[d_1,a_1]$ is connected to
$[d_2,c_2]$. Let $e_1$ be the nearest to $a_1$ vertex in $[d_1,a_1]$
that sends an edge to $[d_2,c_2]$. Let $e_2$ be the nearest to $c_2$ vertex
in $[d_2,c_2]$ that sends an edge to $[d_1,a_1]$. It is not hard to see
that $e_1$ and $e_2$ are adjacent.
Assume that $[d_1,a_1]$ is connected by an edge to $[b_2,c_1]$ and
let $f$ be the nearest to $c_1$ vertex in $[b_2,c_1]$ sending an edge
to $[d_1,a_1]$.

In addition to $s_1$ and $s_2$, remove from $G$ also $e_1$, $e_2$, and $f$,
thereby extending the separator to $\{s_1,s_2,e_1,e_2,f\}$. Then
an inner edge (i.e.\ an edge not in $C$) may be only between
$[d_1,e_1)$ and $[d_2,e_2)$, $(e_1,a_1]$ and $[b_2,f)$,
$(f,c_1]$ and $(e_2,c_2]$. It is easy to see that, if
$[d_1,e_1)$ and $[d_2,e_2)$ are connected by an edge, then
$\flap{d_1}=\flap{d_2}$ is a \eedhop\/ or \edhop\/ flap; otherwise
$\flap{d_1}$ and $\flap{d_2}$ are two flaps, each \hop\/ or \edhop\/
(one of them may disappear if $d_1=e_1$ or $d_2=e_2$).
A completely similar situation is with the pairs $(e_1,a_1]$, $[b_2,f)$
and $(f,c_1]$, $(e_2,c_2]$.

If $[d_1,a_1]$ and $[b_2,c_1]$ are not connected by an edge, we add to
the separator only $e_2$ and, similarly to the above, deal with two
pairs $[d_1,a_1]$, $[d_2,e_2)$ and $(e_2,c_2]$, $[b_2,c_1]$.

In the worst case, the separator we have constructed consists of 5 vertices
and has 7 flaps.
\end{proof}

\section{Bounded degree graphs with separators}\label{s:minor}

\begin{definition}\label{def:simflaps}
Let $G$ be a graph and $X\subset V(G)$. We call two $X$-flaps, $F_1$
and $F_2$, {\em similar\/} if the identity map of $X$ onto itself
extends to an isomorphism from $G[X\cup V(F_1)]$ to $G[X\cup V(F_2)]$.
We say that a separator $\Sep$ for a class of graphs $\cl$ has
{\em at most $s$ similar flaps\/} if, for every $G\in\cl$, the maximum
number of pairwise similar $\sep G$-flaps is bounded by~$s$.
\end{definition}

An $m$-flap separator obviously has at most $m$ similar flaps, and we
hence can expect more classes to have separators with bounded number of
similar flaps rather than with bounded number of all flaps.
We now modify the strategy $S_i$ so that it will guarantee a fast
Spoiler's win on graphs possessing separators with bounded number of
similar flaps. The price will be an increase of the number of alternation
between the graphs during the game.

We denote the modified strategy by $S^*_i$. Like $S_i$, this is a strategy
for Spoiler in the \EF\/ game on non-isomorphic colored graphs $G$
and $G'$, where $G$ is a connected graph in $\col\cl$ and a class of
graphs $\cl$ has separator
$\Sep$. The only difference of $S^*_i$ from $S_i$ concerns Case 2.
In this case Spoiler chooses a flap $H$ such that $m(H)<m'(H)$ and
during $m(H)$
rounds selects one vertex in each $X$-flap isomorphic to $H$.
Recall that all flaps are considered with the modified coloring.
Note that two $X$-flaps with modified colorings can be isomorphic
only if their underlying graphs are similar in the sense of
Definition \ref{def:simflaps}. Thus, $m(H)\le s$, where $s$ is
the maximum number of similar $X$-flaps.

In the next move Spoiler selects a vertex $u$ in $F'$, an $X'$-flap
isomorphic to $H$
that has not been visited by Duplicator, so that $d(u,X')=1$.
If Duplicator responds with a vertex in an $X$-flap isomorphic to $H$,
then, as this $X$-flap already contains a selected vertex, Spoiler
is able to apply the strategy of Lemma \ref{lem:conn} and win in at
most $\lceil\log|H|\rceil$ subsequent moves.
If Duplicator responds with a vertex in an $X$-flap $F$ non-isomorphic
to $H$, then Spoiler restricts the game to $F$ and $F'$ and applies
the strategy $S^*_{i-1}$. If in any round
Spoiler selects a vertex in $F$ but Duplicator replies with a vertex
outside $F'$, Spoiler switches to the strategy of Lemma \ref{lem:conn}
and wins in at most $\lceil\log|F|\rceil$ moves. If
Spoiler selects a vertex in $F'$ but Duplicator replies with a vertex
outside $F$, this is an immediate loss of Duplicator: Lemma \ref{lem:out}
holds true for $S^*_i$ with literally the same proof.

\begin{lemma}\label{lem:length2}
Suppose that a class of connected graphs $\cl$ has hereditary
$\eps$-separator $\Sep$ of size $k(n)$ with at most $s$ similar flaps,
where the function $k$ is
assumed to be defined over reals and non-decreasing.
For $i\ge 0$, let the strategies $S^*_i$ be based on $\Sep$.
If $G$ and $G'$ are colored graphs such that $G\in\col\cl$ and
$G'\not\cong G$, let $L_i(G,G')$ denote the minimum $r$ such that
Spoiler wins $\game_r(G,G')$ following $S_i$, regardless of Duplicator's
strategy. Set $t=\lceil(\log\frac n{s+1})/\log(\eps^{-1})\rceil$. Then
$$
L_t(G,G')<\sum_{i=0}^{t-1}k(\eps^i n)+(s+1)(t+1)+\log n+2.
$$
\end{lemma}

The proof is almost identical to the proof of Lemma~\ref{lem:length}.

\begin{lemma}\label{lem:square}
Let $s\ge 2$, $c>0$, $\eps,\gd\in(0,1)$ be constants such that $\eps(s+1)>1$.
Suppose that $\cl$ is a class of connected graphs with hereditary
$\eps$-separator having at most $s$ similar flaps and size $k(n)=cn^\gd$.
Set $a=2\log n/\log(\eps^{-1})+1$. Then, for all $G\in\cl$, it holds
\begin{equation}\label{eq:square}
\DD aG<
\frac c{1-\eps^\gd}n^\gd + \of{\frac{s+1}{\log(\eps^{-1})}+1}\log n+s+3.
\end{equation}
\end{lemma}

\begin{proof}
It is easy to see that, following $S^*_i$, Spoiler during the game makes
at most $2i+1$ alternations between $G$ and $G'$. Similarly to the
proof of Lemma \ref{lem:const}, we conclude from Lemma \ref{lem:length2}
that
$$
\DD{2t+1}G<
\sum_{i=0}^{t-1}k(\eps^in)+(s+1)t+\log n+s+3,
$$
where
$$
t<\frac{\log\frac n{s+1}}{\log(\eps^{-1})}+1\le\frac{\log n}{\log(\eps^{-1})}
$$
(the latter inequality by the condition $\eps(s+1)>1$).
To obtain \refeq{eq:square}, it remains to notice that $2t+1<a$ and
$\sum_{i=0}^{t-1}k(\eps^in)=cn^\gd(1-\eps^{\gd t})/(1-\eps^\gd)$.
\end{proof}

Seeking for applications of Lemma \ref{lem:square}, we consider
the following classes of graphs for which $\eps$-separators
of size $O(\sqrt n)$ are known.

{\it Planar graphs}.
An efficient construction of $\frac23$-separators of size $2\sqrt2\sqrt n$
is given in \cite{LTa}. The constant $2\sqrt2$ is improved to $\frac32\sqrt2$
in~\cite{AST0}.

{\it Graphs of genus $g$ (embeddable in a surface of genus $g$)}.
An efficient construction of $\frac23$-separators of size $O(\sqrt{gn})$
is given in~\cite{GHT}.

{\it Graphs with no $K_h$ minor}.
An efficient construction of $\frac12$-separators of size $h^{3/2}\sqrt{n}$
is given in~\cite{AST}.

To see inclusions between these classes, note that planar graphs are,
in other words, graphs of genus 0, and that a graph of genus $g$
cannot have a $K_h$ minor if $g<(h-3)(h-4)/12$.
The Robertson-Seymour Graph Minor Theorem implies
that every proper class of graphs closed under minors is, for some $h$,
contained in the class of graphs with excluded $K_h$ minor.

Unfortunately, for the above classes no bound $o(n)$ for the number of
similar flaps is possible. We therefore need to put an additional
restriction.

\begin{lemma}\label{lem:simflaps}
Let $\Sep$ be a separator for a class of graphs $\cl$.
Let $\cl_d$ denote the restriction of $\cl$ to connected graphs
of degree at most $d$. Then $\Sep$ is a separator for $\cl_d$
with at most $d$ similar flaps.
\end{lemma}

\begin{proof}
Suppose that $G\in\cl_d$ and $X=\sep G$. Since $G$ is connected,
all $X$-flaps are connected by an edge to $X$.  Pairwise similar $X$-flaps,
by definition, are all connected by an edge to the same vertex of $X$.
The maximum number of such flaps is therefore at most~$d$.
\end{proof}

\begin{theorem}\label{thm:minor}
Let $H(G)$ denote the smallest $h$ such that a graph $G$ does not
have a minor $K_h$. Then, for every connected $G$,
\begin{equation}\label{eq:all}
\DD{2\log n+1}G<
(2+\sqrt2)H(G)^{3/2}\sqrt n+(\Delta(G)+2)(\log n+1)+1.
\end{equation}
In particular, if $H(G)=O(1)$ and $\Delta(G)=O(\sqrt n/\log n)$, then
$$
\DD{2\log n+1}G=O(\sqrt n).
$$
If $G$ is connected planar, then
\begin{equation}\label{eq:planar}
\DD{2\log n/\log(3/2)+1}G<
(\frac92\sqrt2+3\sqrt3)\sqrt n+\of{\frac{\Delta(G)+1}{\log(3/2)}+1}\log n
+\Delta(G)+3.
\end{equation}
If $G$ is connected and has genus $g$, then
\begin{equation}\label{eq:genus}
\DD{2\log n/\log(3/2)+1}G<
c\sqrt g\sqrt n+\of{\frac{\Delta(G)+1}{\log(3/2)}+1}\log n
+\Delta(G)+3,
\end{equation}
where $c$ is a constant.
\end{theorem}

\begin{proof}
Let us prove the bound \refeq{eq:all}. Fix an arbitrary connected graph
$G$ and denote $h=H(G)$ and $d=\Delta(G)$. Let $\cl(G)$ consists of
all connected induced subgraphs of $G$. Clearly, every graph in $\cl(G)$
has degree at most $d$ and no minor $K_h$. By \cite{AST} and Lemma
\ref{lem:simflaps}, graphs in $\cl(G)$ have $\frac12$-separators
of size $h^{3/2}\sqrt n$ with at most $d$ similar flaps. Since $\cl(G)$
is closed under connected induced subgraphs, the separator is hereditary.
Thus, Lemma \ref{lem:square} applies and completes the proof.
The bounds \refeq{eq:planar} and \refeq{eq:genus} follow in the same
way from the separator theorems of \cite{AST0} and~\cite{GHT}.
\end{proof}

\end{document}